\documentclass{amsart}
\usepackage[dvips]{graphicx}
\usepackage{amscd}
\usepackage{amsmath}
\usepackage{amsfonts}
\usepackage{amssymb}
\newtheorem{theorem}{Theorem}
\theoremstyle{plain}
\newtheorem{acknowledgement}{Acknowledgement}

\newtheorem{lemma}{Lemma}

\newtheorem{problem}{Problem}

\numberwithin{equation}{section}

\begin{document}
\title{Ringel's generalized earth-moon problem}
\author{Jeremy Alexander}
\address{Mathematics Department SCEN 301\\
University of Arkansas\\
Fayetteville, AR 72701}
\email{jralexa@uark.edu}
\urladdr{comp.uark.edu/\symbol{126}jralexa}
\date{April 19, 2005}
\maketitle

\section{\bigskip Introduction}

A graph G has surface thickness n (S-thickness n) with respect to the surface
S, if G can be decomposed into n and no fewer subgraphs by making n copies of
the vertex set of G and then assigning each edge of G to one of the n copies
so that n graphs result; each resulting subgraph of G must be S embeddable.
\ The chromatic number of a graph G, denoted $\chi$(G), is simply the fewest
number of colors needed to color the vertex set of G such that no two adjacent
vertices receive the same color. \ Define $\chi_{\text{\textit{n}}}$(S) to be
the chromatic number of the surface S, that is, the minimum number of colors
needed to properly color all S-thickness\ n graphs. \ Let E(S) denote the
Euler characteristic of the surface S. \ Let S$_{\text{k}}$ denote the
orientable surface of genus k and N$_{\text{k}}$ the nonorientable surface
obtained by adding k crosscaps to the sphere; S$_{\text{0}}$ denotes the
two-dimensional sphere. \ It is known that E(S$_{\text{k}}$) = 2 - 2k and
E(N$_{\text{k}}$) = 2 - k.

\bigskip

In 1959 [1, p.233], Ringel asked: \ What is the chromatic number of the sphere
for graphs of thickness two? \ The bounds are known to be: \ 9 $\leq$
$\chi_{\text{\textit{2}}}$(S$_{\text{0}}$) $\leq$ 12. \ See [2] for a
delightful exposition of the spherical case. \ In [1], Ringel and Jackson ask:
\ What is $\chi_{\text{\textit{n}}}$(S) for any surface S? \ There [p.241] the
following upper bound is given for any surface S, except the sphere:

\bigskip

\begin{center}
$\chi_{\text{\textit{n}}}$(S) $\leq$ $\left\lfloor \text{( 1 + 6n +}%
\sqrt{\text{ ( 1 + 6n )}^{2}\text{ - \ 24nE(S)}}\text{ ) / 2}\right\rfloor $ \ \ \ \ \ \ \ \ \ \ \ \ \ \ \ (1)

\bigskip
\end{center}

The main results in this note are three elementary arguments that establish
$\chi_{\text{\textit{n}}}$(S) for two new surfaces: \ N$_{\text{3}}$, also
known as Dyck's surface, and the Moebius band which will be denoted by M.
\ The arguments are based on previous results appearing in [1]. \ The results
are not especially surprising, but they do provide good lower bounds for
$\chi_{\text{\textit{n}}}$(N$_{\text{2k+1}}$) once the value of $\chi
_{\text{\textit{n}}}$(S$_{\text{k}}$) is known for k
$>$%
0; in the case k = 1, the exact value of $\chi_{\text{\textit{n}}}%
$(N$_{\text{3}}$) is obtained.

\pagebreak 

\section{Main Results}

\begin{lemma}
$\left\lfloor \text{( 1+6n +}\sqrt{\text{ ( 1+6n )}^{2}\text{ - \ 24nE(N}%
_{\text{3}}\text{)}}\text{ ) / 2}\right\rfloor $ $\leq$ 6n+1 for all n $\geq$ 2.
\end{lemma}

\begin{proof}
Suppose for contradiction that the claim is false.

Then since E(N$_{\text{3}}$) = -1, there exists an integer n $\geq$ 2 such that:

$\left\lfloor \text{( 1 + 6n +}\sqrt{\text{ ( 1 + 6n )}^{2}\text{ + \ 24n}%
}\text{ ) / 2}\right\rfloor $
$>$%
6n + 1 $\rightarrow$

$\left\lfloor \text{( 1 + 6n +}\sqrt{\text{ ( 1 + 6n )}^{2}\text{ + \ 24n}%
}\text{ ) / 2}\right\rfloor $ $\geq$ 6n + 2 $\rightarrow$

1 + 6n +$\sqrt{\text{ ( 1 + 6n )}^{2}\text{ + \ 24n}}$ $\geq$ 12n + 4
$\rightarrow$

$\sqrt{\text{ ( 1 + 6n )}^{2}\text{ + \ 24n}}$ $\geq$ 6n + 3 $\rightarrow$

( 1 + 6n )$^{2}$ + \ 24n $\geq$ ( 6n + 3 )$^{2}$ $\rightarrow$

1 $\geq$ 9. \ This is an obvious contradiction, hence, the original assumption

was false and the lemma is true.
\end{proof}

\begin{theorem}
\bigskip$\chi_{\text{\textit{n}}}$(N$_{\text{3}}$) = 6n + 1 for all n $\geq$ 2.
\end{theorem}

\begin{proof}
In [1, p.238], it is shown that $\chi_{\text{\textit{n}}}$(S$_{\text{1}}$) =
6n + 1 for all n $\geq$ 2. \ Let G be any graph embedded on S$_{\text{1}}$.
\ Choose any face of S$_{\text{1}}$ created by the embedding of G and add a
crosscap, what can be obtained is an embedding of G on the surface
N$_{\text{3}}$ since N$_{\text{3}}$ is homeomorphic to the surface
S$_{\text{1}}$ + crosscap and homeomorphism preserves graph isomorphism as
well as graph embedding. \ Combining this with the previous lemma gives the
desired result.
\end{proof}

For a complete classification of surfaces see Conway's ZIP proof in [3].

\begin{theorem}
$\chi_{\text{\textit{n}}}$(M) = 6n for all n $\geq$ 2.
\end{theorem}

\begin{proof}
M is homeomorphic to N$_{\text{1}}$ minus an open disk, to see this look at
the figure below; a square minus the region X (an open disk) with opposite
ends identified. \ Notice that coupled ends have opposite orientation and so a
half-twist is needed to join each pair correctly. \ When opposite ends are
joined a Moebius band results. \ Now, fill in the two missing regions denoted
by X. \ If sides of the resulting square receive the same identifications
induced by the labeling of the figure below, N$_{\text{1}}$, the projective
plane, is obtained. \ Now, given any graph G embedded on N$_{\text{1}}$,
select any face induced by the embedding and delete it, what results is an
embedding of G on a Moebius band. \ Inductively it follows that $\chi
_{\text{\textit{n}}}$(M) $\geq$ 6n for all n $\geq$ 2 since in [1, p.240] it
is proven that $\chi_{\text{\textit{n}}}$(N$_{\text{1}}$) = 6n for all n
$\geq$ 2. \ Now, inequality (1) gives an upper bound of $\chi
_{\text{\textit{n}}}$(M) $\leq$ 6n + 1 for all n $\geq$ 2 since E(M) = 0. \ To
show that this upper bound cannot be obtained, the case $\chi
_{\text{\textit{2}}}$(M) is shown, from which the general result easily
follows. \ Suppose for contradiction that $\chi_{\text{\textit{2}}}$(M) = 13.
\ Then there exists a graph G of M-thickness two such that\ $\chi$(G) = 13.
\ Since the square minus an open disk appearing in the picture below is
homeomorphic to M, G biembeds on this version of a Moebius band. \ Now fill in
the deleted disk labeled X in the figure below to obtain a biembedding of the
graph G on N$_{\text{1}}$. \ This implies that $\chi_{\text{\textit{2}}}%
$(N$_{\text{1}}$) = 13 which is a contradiction since $\chi_{\text{\textit{2}%
}}$(N$_{\text{1}}$) = 12.
\end{proof}

\bigskip\pagebreak 

\bigskip%
\begin{center}
\includegraphics[
height=3.6876in,
width=3.6685in
]%
{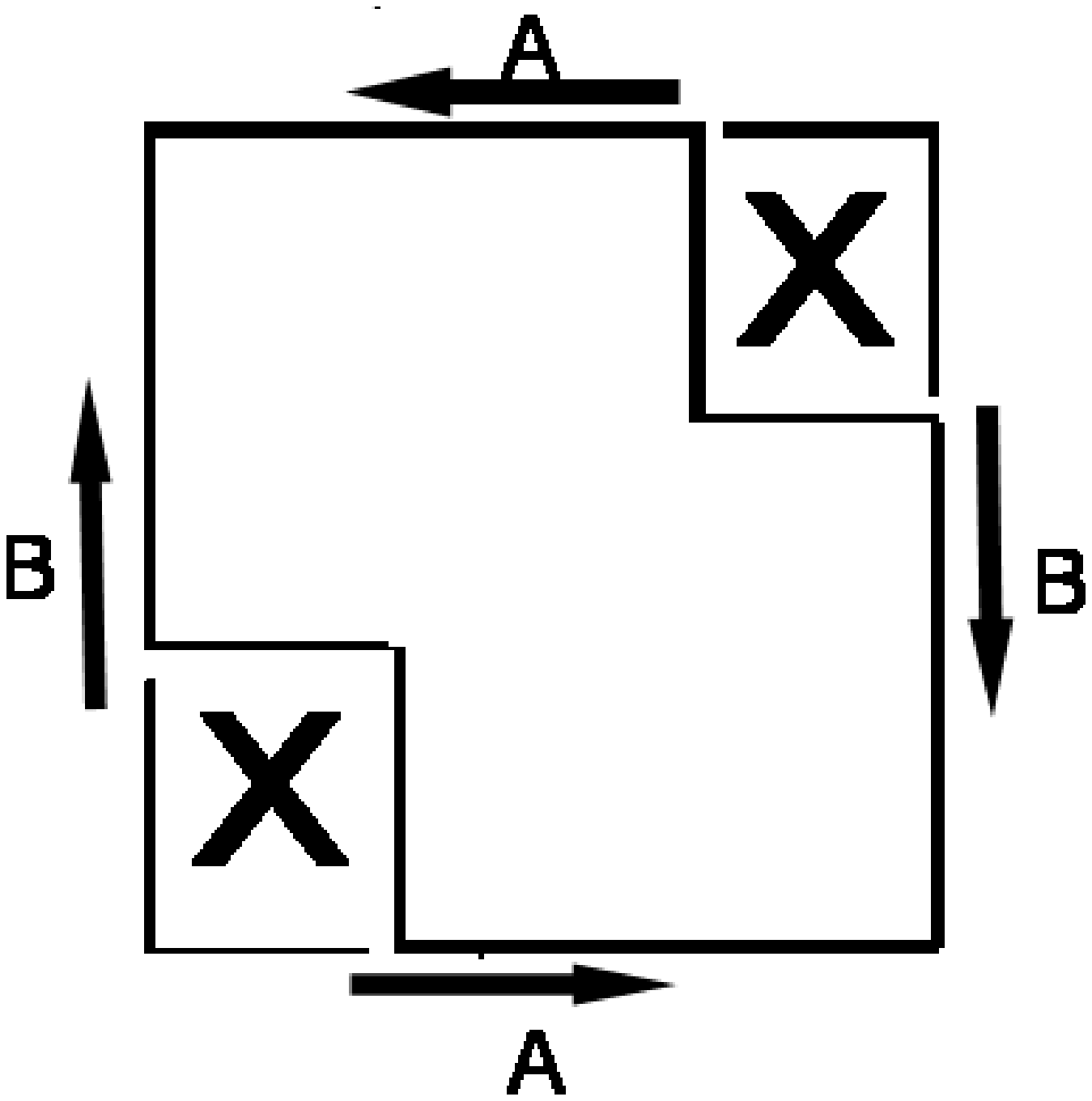}%
\end{center}

Hence, Ringel's generalized earth-moon problem has been solved for two new
surfaces M and N$_{\text{3}}$, also known as the Moebius band and Dyck's
surface, respectively. \ The following list summarizes known results for
Ringel's generalized earth-moon problem:

For all n $\geq$ 2,

$\chi_{\text{\textit{n}}}$(S$_{\text{2}}$) = 6n + 2 \ (see [5])

$\chi_{\text{\textit{n}}}$(S$_{\text{1}}$) = 6n + 1 \ (see [1], [5])

$\chi_{\text{\textit{n}}}$(N$_{\text{1}}$) = 6n \ (see [1], [5])

$\chi_{\text{\textit{n}}}$(M) = 6n

$\chi_{\text{\textit{n}}}$(N$_{\text{3}}$) = 6n + 1.

Until very recently the N$_{\text{2}}$-thickness of K$_{\text{13}}$ was not
known, however, it appears that Thom Sulanke has determined it is three using
an exhaustive computer search [4]. \ The N$_{\text{2}}$-thickness of
K$_{\text{13}}$ has been a long standing obstruction to extending Beineke's
results appearing in [5] for the Klein Bottle thickness of K$_{\text{n}}$. \ A
similar obstruction for determining the S$_{\text{3}}$-thickness of
K$_{\text{n}}$ is the S$_{\text{3}}$-thickness of K$_{\text{16}}$ [5, p.996];
Sulanke has determined that $\chi_{\text{\textit{2}}}$(S$_{\text{3}}$) = 16,
although this result has not been published. \ Through personal communication
Thom has also told me that he has been able to determine $\chi
_{\text{\textit{n}}}$(N$_{\text{4}}$) and $\chi_{\text{\textit{n}}}%
$(N$_{\text{5}}$). \ He also pointed out to me that Lemma 1, appearing here,
generalizes for all integers n
$>$%
0 and E(S)
$<$%
0. \ Using the exact same arguments given in the proof of Lemma 1, one can
show that the left side of inequality (1) is $\leq$ 6n - E(S) for n
$>$%
0 and E(S)
$<$%
0. \ Naturally, the next two surfaces to consider are S$_{\text{3}}$ and
N$_{\text{2}}$. \ This motivates the following problem:

\begin{problem}
Determine $\chi_{\text{\textit{n}}}$(S$_{\text{3}}$) and $\chi
_{\text{\textit{n}}}$(N$_{\text{2}}$).
\end{problem}

\begin{acknowledgement}
I would like to thank Yo'av Rieck and Chaim Goodman-Strauss for constructive criticism.
\end{acknowledgement}

\pagebreak 

\section{References}

[1] B. Jackson, G. Ringel, Variations on Ringel's earth-moon problem,

\textit{Discrete Mathematics} \textbf{211} (2000), 233 - 242.

[2] J. P. Hutchinson, Coloring ordinary maps, maps of empires, and maps

of the moon, \textit{Math. Magazine} (4) \textbf{66} (1993), 211 - 226.

[3] G. Francis, J. Weeks, Conway's ZIP proof,

\textit{Amer. Math. Monthly} \textbf{106} (1999), 393-399.

[4] T. Sulanke, Biembeddings of K$_{\text{13}}$,

\textit{http://needmore.physics.indiana.edu/\symbol{126}tsulanke/ graphs/biembed/biembed.pdf}

[5] L. W. Beineke, Minimal decompositions of complete graphs into

subgraphs with embeddability properties, \textit{Canad. J. Math} \textbf{21} (1969),

992-1000.
\end{document}